\documentclass[11pt]{article}
\usepackage{dsfont}
\usepackage{graphics}
\usepackage[demo]{graphicx}
\usepackage{epsfig}
\usepackage{pstricks}
\usepackage[normal]{subfigure}
\usepackage[latin1]{inputenc}
\usepackage[english,francais]{babel}
\usepackage{relsize,exscale}
\usepackage{makeidx}
\usepackage{enumitem}
\usepackage{amsfonts,amssymb,amsmath}
\usepackage{graphicx}
\usepackage{color}
\usepackage{multirow}
\usepackage{mathrsfs}
\usepackage[normalem]{ulem}
\usepackage{cancel}
\usepackage{bbm}
\newenvironment{prooff}{{\it Proof :}}{\hfill\rule{2mm}{2mm}\vskip3mm\par}
\newtheorem{theorem}{Theorem}[section]

\newtheorem{corollary}[theorem]{Corollary}
\newtheorem{e-definition}[theorem]{Definition\rm}
\newtheorem{remark}{\it Remark\/}

\usepackage{color}
\definecolor{dred}{rgb}{0.92,0,0}
\definecolor{dgreen}{rgb}{0,0.92,0}
\definecolor{dblue}{rgb}{0,0,0.92}
\definecolor{dyellow}{rgb}{0.95,0.95,0}

\newcommand{\R}{\mathbb{R}}
\newcommand{\N}{\mathbb{N}}
\def\D{\displaystyle}
\newcommand{\hs}{\hspace{0.1cm}}
\newcommand{\vs}{\vspace{0.4cm}}
\newcommand{\sa}{\\ [0.2cm]}
\usepackage{multirow}
\graphicspath{
{./Figures/}
{./}
}
\title{A probabilistic approach for exact solutions of determinist PDE's as well as their finite element approximations}
\author{Jo\"el Chaskalovic \\[0.2cm]
D'Alembert,
Sorbonne University, Paris, France, \\[0.1cm]\emph{email}: jch1826@gmail.com
}
\date{}
\begin{document}
\maketitle
\selectlanguage{english}
\begin{abstract}
\noindent A probabilistic approach is developed for the exact solution $u$ to a determinist partial differential equation as well as for its associated approximation $u^{(k)}_{h}$ performed by $P_k$ Lagrange finite element. Two limitations motivated our approach: on the one hand, the inability to determine the exact solution $u$ to a given partial differential equation (which initially motivates one to approximating it) and, on the other hand, the existence of uncertainties associated with the numerical approximation $u^{(k)}_{h}$. We thus fill this knowledge gap by considering the exact solution $u$ together with its corresponding approximation $u^{(k)}_{h}$ as random variables. By way of consequence, any function where $u$ and $u_{h}^{(k)}$ are involved as well. In this paper, we focus our analysis to a variational formulation defined on $W^{m,p}$ Sobolev spaces and the corresponding \emph{a priori} estimates of the exact solution $u$ and its approximation $u^{(k)}_{h}$ to consider their respective $W^{m,p}-$norm as a random variable, as well as  the $W^{m,p}$ approximation error with regards to $P_k$ finite elements. This will enable us to derive a new probability distribution to evaluate the relative accuracy between two Lagrange finite elements $P_{k_1}$ and $P_{k_2}, (k_1 < k_2)$.\sa
\noindent{\footnotesize {\em keywords}: Error estimates, Finite elements, Bramble-Hilbert lemma, Sobolev spaces.}
\end{abstract}
\section{Introduction}\label{intro}
\noindent We recently proposed new perspectives on relative finite element accuracy (\cite{ChasAs18} and \cite{Arxiv2}), using a mixed geometrical-probabilistic interpretation of the error estimate in the case of finite element approximation (see for example \cite{Ciarlet} or \cite{RaTho82}), derived from Bramble-Hilbert lemma \cite{Bramble_Hilbert}. \sa
In \cite{arXiv_Wmp}, we further extended the results we had derived in the case of $H^1$ to Sobolev spaces $W^{m,p}, (p\neq 2)$. To this end, we had to consider a more general framework, which mainly relied on the Banach-Ne$\check{\mbox{c}}$as-Babus$\check{\mbox{h}}$ka $\textbf{(BNB)}$ abstract problem \cite{Ern_Guermond} devoted to Banach spaces. \sa
This enabled us to obtain two new probability distributions which estimate the relative accuracy between two Lagrange finite elements $P_{k_1}$ and $P_{k_2}$, ($k_1 < k_2$), by considering it as a random variable.\sa
We thus obtained new results which show, amongst others, which of $P_{k_1}$ or $P_{k_2}$ is the most likely accurate, depending on the value of the mesh size $h$; this value is not considered anymore as going to zero, as in the standard procedure.\sa
However, while obtaining these probability distributions, we only considered the standard error estimate dedicated to $P_k$ Lagrange finite elements, which approximates the solution $u$ to a variational problem $\textbf{(BNB)}$, $u$ being formulated in the present case in the Sobolev space $W^{m,p}(\Omega)$. \sa
In the current work, we enrich the model published in \cite{arXiv_Wmp} in multiple manners. Indeed, considering the functional framework $\textbf{(BNB)}$ in the case of the $W^{m,p}$ Sobolev spaces, we will take into account the available \emph{a priori} estimates one can deal with the solution of this kind of problem as well as for its approximation, together with the approximation error which corresponds to $P_k$ Lagrange finite elements (see Section \ref{Wmp_estimate}). This will enable us to derive a new probabilistic model applied to the relative accuracy between two finite elements $P_{k_1}$ and $P_{k_2}, (k_1<k_2)$. To this end, we also generalize the discrete probabilistic framework we considered in \cite{ChasAs18} and \cite{arXiv_Wmp} by introducing a continuous probabilistic formalism based on an appropriate density of probability (see Section \ref{FEM_accuracy}).\sa
The paper is organized as follows. In Section \ref{Second_Order_Elliptic}, we recall the mathematical problem we consider and introduce the basic definitions of functional tools to consider different estimations in $W^{m,p}$ Sobolev spaces. Section \ref{FEM_accuracy} is dedicated to the analysis of the relative finite elements accuracy based on a probabilistic approach. In Section \ref{discussion}, we detail the contribution of the \emph{a priori} estimates and their interactions with the error estimate in the probability distributions we derived in Section \ref{FEM_accuracy}. Concluding remarks follow.
\section{The functional and approximation frameworks and their corresponding estimates}\label{Second_Order_Elliptic}
\subsection{Abstract problem in Banach spaces and corresponding fundamental results}
\noindent In this section we define an abstract framework which will enable us to consider the solution of a variational problem in Banach $W^{m,p}$ Sobolev spaces when $p\neq 2$, and its corresponding approximation computed by $P_k$ Lagrange finite elements.\sa
\noindent In order to do so, we follow the presentation of A. Ern and J. L. Guermond \cite{Ern_Guermond}, where two general Banach spaces $W$ and $V$ are involved with $V$ reflexive. We also recall the different assumptions needed in order to apply the $\textbf{(BNB)}$ Theorem {\sout valid} in Banach spaces.\sa
Let $u \in W$ be the solution of the variational formulation $\textbf{(VP)}$ defined by:
\begin{equation}\label{VP}
\textbf{(VP}\textbf{)} \hspace{0.2cm} \left\{
\begin{array}{l}
\mbox{Find } u \in   W \mbox{ such that:} \\ [0.1cm]
a(u,v) = l(v), \quad\forall v \in V, 
\end{array}
\right.
\end{equation}
where:
\begin{enumerate}
\item $W$ and $V$ are two Banach spaces equipped with norms denoted by $\|.\|_W$ and $\|.\|_V$, respectively; moreover, $V$ is reflexive.\label{Banach}
\item $a$ is a continuous bilinear form on $W\times V$, i.e, $a\in\mathcal{L}(W\times V;\R):$ $$\forall (u,v)\in W\times V,\, |a(u,v)|\leq\|a\|_{W,V}\|u\|_W\|v\|_V,$$
with: $\D\|a\|_{W,V}\equiv\inf\left\{C\in\R^{*}_{+},\forall (u,v)\in W\times V: |a(u,v)|\leq C\|u\|_W\|v\|_V\right\}$.\label{Norme_a}
\item $l$ is a continuous linear form on $V$, i.e, $l\in V'=\mathcal{L}(V;\R)$.
\end{enumerate}
We further make the two following assumptions:
\begin{description}
\item[\textbf{(BNB1)}] \hspace{4cm}$\D \exists \alpha > 0, \hs \inf_{w\in W}\sup_{v\in V}\frac{a(w,v)}{\|w\|_{W}\|v\|_{V}} \geq \alpha $,\vspace{0.2cm}
\item[\textbf{(BNB2)}] \hspace{3cm}$\forall v\in V, (\forall w\in W, a(w,v)=0)\Longrightarrow (v=0)$.
\end{description}
\vs
Then, one can prove the \textbf{(BNB)} Theorem (\cite{Ern_Guermond}, Theorem 2.6) which claims that variational problem \textbf{(VP)} has one and only one solution in $W$ and that the following \emph{a priori} estimate holds:
\begin{equation}\label{estimation_a_priori}
\forall \l\in V': \|u\|_{W} \leq \frac{\|l\|_{V'}}{\alpha}.
\end{equation}
We also define the approximation $u^{(k)}_{h}$ of $u$, solution to the approximate variational formulation:
\begin{equation}\label{VP_h}
\textbf{(VP}\textbf{)}_{h} \hspace{0.2cm} \left\{
\begin{array}{l}
\mbox{Find } u^{(k)}_{h} \in   W^{(k)}_h \mbox{ such that:} \\ [0.15cm]
a(u^{(k)}_{h},v^{(k)}_{h}) = l(v^{(k)}_{h}),\quad \forall v^{(k)}_{h} \in V^{(k)}_h, 
\end{array}
\right.
\end{equation}
where we assume that $W^{(k)}_h$ and $V^{(k)}_h$ are two finite-dimensional subsets of $W$ and $V$, respectively.\sa
Moreover, as noticed in \cite{Ern_Guermond} (Remark 2.23, p.92), neither condition \textbf{(BNB1)} nor condition \textbf{(BNB2)} imply its discrete counterpart. Then, the well-posedness of (\ref{VP_h}) is equivalent to the two following discrete conditions:
\begin{description}
\item[(BNB1$_h$)] \hspace{4cm}$\D \exists \alpha^{(k)}_h > 0, \hs \inf_{w^{(k)}_h\in\, W^{(k)}_h}\sup_{v^{(k)}_h\in\, V^{(k)}_h}\frac{a(w^{(k)}_h,v^{(k)}_h)}{\|w^{(k)}_h\|_{W^{(k)}_h}\|v^{(k)}_h\|_{V^{(k)}_h}} \geq \alpha^{(k)}_h$,\vspace{0.2cm}
\item[(BNB2$_h$)] \hspace{3cm}$\forall v^{(k)}_h\in V^{(k)}_h, (\forall w^{(k)}_h\in W^{(k)}_h, a(w^{(k)}_h,v^{(k)}_h)=0)\Longrightarrow (v^{(k)}_h=0)$.
\end{description}
\vspace{0.2cm}
If we furthermore assume that dim$(W^{(k)}_h)$ = dim$(V^{(k)}_h)$, a direct application of Theorem 2.2 in \cite{Ern_Guermond} enables us to write the following \emph{a priori} estimate:
\begin{equation}\label{estimation_a_priori_discret}
\forall \l\in V': \|u_{h}^{(k)}\|_{W_h} \leq \frac{\|l\|_{V'}}{\alpha^{(k)}_h}.
\end{equation}
In the next section we will apply these results to the particular case where the exact solution $u$ belongs to $W^{m,p}$ Sobolev spaces and the approximation $u^{(k)}_{h}$ is computed by the help of $P_k$ Lagrange finite elements.
\subsection{Application to $W^{m,p}$ Sobolev spaces and the corresponding error estimate}\label{Wmp_estimate}
\noindent We introduce an open-bounded subset $\Omega \subset \R^n$ exactly recovered by a mesh ${\mathcal T}_h$ composed by $N_K$ $n$-simplexes $K_{\mu}, (1\leq \mu \leq N_K),$ which respect the classical rules of regular discretization (see for example \cite{RaTho82}). We moreover denote by $h$ the mesh size of ${\mathcal T}_h$ (the largest diameter in the mesh ${\mathcal T}_h$), and by $P_k(K_{\mu})$ the space of polynomials defined on a given $n$-simplex $K_{\mu}$ of degree less than or equal to $k$ ($k \geq$~1). \sa
Thenceforth, we assume that the approximate spaces $W^{(k)}_h$ and $V^{(k)}_h$ satisfy dim$(W^{(k)}_h)$ = dim$(V^{(k)}_h)$, that they are included in the space of functions defined on $\Omega$, and composed of polynomials belonging to $P_k(K_{\mu}), (1 \leq \mu \leq N_K)$. \sa
Finally, we also specify the functional framework of the abstract problem \textbf{(VP)} by introducing $W^{m,p}$ Sobolev spaces as follows: \sa
\noindent For any integer $m\geq 1$ and any $1\leq p \leq +\infty$, we denote by $W^{m,p}(\Omega)$ the Sobolev space of (class of) real-valued functions which, together with all their partial distributional derivatives of order less or equal to $m$, belongs to $L^p(\Omega)$:
\begin{equation}
\D W^{m,p}(\Omega) = \left\{\!\!\frac{}{}u \in L^p(\Omega)\,/\,\forall\, \alpha, |\alpha|\leq m, \partial^{\alpha}u\in L^p(\Omega)\right\},
\end{equation}
$\alpha=(\alpha_1, \alpha_2, \ldots, \alpha_n) \in \N^{n}$ being a multi-index whose length $|\alpha|$ is given by $|\alpha|=\alpha_1+\dots+\alpha_n$, and $\partial^{\alpha}u$ being the partial derivative of order $|\alpha|$ defined by:
\begin{equation}
\D \partial^{\alpha}u \equiv \frac{\partial^{|\alpha|}u}{\partial x_{1}^{\alpha_1}\dots\partial x_{n}^{\alpha_n}}.
\end{equation}
We also consider the norm $\|.\|_{m,p,\Omega}$ and the semi-norms $|.|_{l,p,\Omega}$, which are respectively defined by:
\begin{equation}
\D \forall u \in\,W^{m,p}(\Omega): \|u\|_{m,p,\Omega} = \left(\sum_{|\alpha|\leq m}\|\partial^{\alpha}u\|^{p}_{L^p}\right)^{1/p}\!\!\!\!, \hs\hs|u|_{l,p,\Omega} = \left(\sum_{|\alpha|= l}\|\partial^{\alpha}u\|^{p}_{L^p}\right)^{1/p}\!\!\!\!, 0 \leq l \leq m,
\end{equation}
where $\|.\|_{L^p}$ denotes the standard norm in $L^p(\Omega)$. \sa
Then, in order to fulfill the conditions of the \textbf{(BNB)} Theorem, particularly so that $W$ be a Banach space and $V$ a reflexive one, in the sequel of the paper, the following definitions of spaces $W$ and $V$ hold:
\begin{equation}\label{V_and_W}
W\equiv W^{m,p}(\Omega) \mbox{ and } V\equiv W^{m',p'}(\Omega),
\end{equation}
where $m$ and $m'$ are two non zero integers and $p$ and $p'$ two real positive numbers which satisfy $p\ne 2$ and $p'>1$ such that:
\begin{equation}\label{Conjugated}
\D \frac{1}{p}+\frac{1}{p'}=1.
\end{equation}
Regarding these choices, Sobolev's space $\left(W^{m,p}(\Omega),\|.\|_{m,p,\Omega}\right)$ is a Banach space and $\left(W^{m',p'}(\Omega),\|.\|_{m',p',\Omega}\right)$ is a reflexive one \cite{Brezis}. Moreover, we have: $W^{(k)}_h \subset W^{m,p}(\Omega)$ and $V^{(k)}_h \subset W^{m',p'}(\Omega)$. \sa
We can now recall the $W^{m,p}$ \emph{a priori} error estimate for $P_k$ Lagrange finite elements we derived in \cite{arXiv_Wmp}:
\begin{equation}\label{estimation_error}
\|u- u^{(k)}_h\|_{m,p,\Omega} \hs \leq \hs \mathscr{C}_k\,h^{k+1-m} \, |u|_{k+1,p,\Omega}\,,
\end{equation}
\vspace{0.1cm}
where $\mathscr{C}_k$ is a positive constant independent of $h$.
\begin{remark}
Since we noticed that $W^{(k)}_h$ is included in $W^{m,p}(\Omega)$, by considering for the topology of $W^{(k)}_h$ that induced from $W^{m,p}(\Omega)$, thanks to the triangle inequality, (\ref{estimation_a_priori})- (\ref{estimation_a_priori_discret}) and (\ref{estimation_error}) lead to the following error estimate:
\begin{equation}\label{Error_Estimate_Approx_0}
\|u- u^{(k)}_h\|_{m,p,\Omega} \hs \leq \hs \min\left(\left[\frac{1}{\alpha}+\frac{1}{\alpha_h^{(k)}}\right]\|l\|_{V'},\,\mathscr{C}_k\,h^{k+1-m} \, |u|_{k+1,p,\Omega}\right),
\end{equation}
where here $V'$ is the dual of $W^{m'\!,p'}(\Omega)$.\sa
As one can see, the right-hand side of (\ref{Error_Estimate_Approx_0}) contains $\alpha_h^{(k)}$, whose dependency on $h$ is usually unknown, except for particular cases. This is the reason why we will assume it is bounded form below by a positive constant $\delta$ independent of h:
$$\forall h >0, 0 < \delta \leq \alpha_h^{(k)}.$$
This uniform boundedness property of $\alpha_h^{(k)}$, crucial to guarantee optimal error estimates \cite{Ern_Guermond}, is valid in multiple cases (see Chapters 4 and 5 in \cite{Ern_Guermond} or \cite{Toundykov_Avalos}), but not systematically, (see for example first-order PDE's in \cite{Ern_Guermond}). \sa
Next, from (\ref{Error_Estimate_Approx_0}) we obtain:
\begin{equation}\label{Error_Estimate_Approx}
\|u- u^{(k)}_h\|_{m,p,\Omega} \hs \leq \hs \min\left(\frac{\|l\|_{V'}}{\alpha^*},\,\mathscr{C}_k\,h^{k+1-m} \, |u|_{k+1,p,\Omega}\right),
\end{equation}
where $\D\frac{1}{\alpha^*}=\frac{1}{\alpha}+\frac{1}{\delta}$.\sa
In the sequel, we will denote by $\beta_k$ the following expression:
\begin{equation}\label{Lambda}
\beta_k\equiv\min\left(\frac{\|l\|_{V'}}{\alpha^*},\,\mathscr{C}_k\,h^{k+1-m} \, |u|_{k+1,p,\Omega}\right).
\end{equation}
Finally, we also observe that there exists a critical value $\hslash_k$ of $h$ defined by:
\begin{equation}\label{hslash_k}
\D \hslash_k \equiv \left(\frac{\|l\|_{V'}}{\alpha^*\,\mathscr{C}_k\,|u|_{k+1,p,\Omega}}\right)^{\frac{1}{k+1-m}},
\end{equation}
such that the error estimate (\ref{Error_Estimate_Approx}) can be splitted according to:
\begin{equation}\label{Splitted_Error_Estimate_Approx}
\|u- u^{(k)}_h\|_{m,p,\Omega}  \leq \left |
\begin{array}{ll}
\D \frac{\|l\|_{V'}}{\alpha^*}\,, & \mbox{ if }\, h \geq \hslash_k, \medskip \\
\D \mathscr{C}_k\,h^{k+1-m} \, |u|_{k+1,p,\Omega}\,, & \mbox{ if }\, h \leq \hslash_k.
\end{array}
\right.
\end{equation}
\end{remark}
Based on this error estimate, the following section is devoted to derive a probabilistic model applied to the relative accuracy between two Lagrange finite elements $P_{k_1}$ and $P_{k_2}, (k_1<k_2),$ when the mesh size $h$ has a given and fixed value.
\section{The probabilistic analysis of the relative finite elements accuracy}
\noindent In \cite{arXiv_Wmp} we proposed two probability distributions which enabled us to appreciate the evaluation of the more likely $W^{m,p}-$accurate between two Lagrange finite elements $P_{k_{1}}$ and $P_{k_{2}}, (k_1 < k_2)$. These distributions were essentially derived by considering the error estimate (\ref{estimation_error}).\sa
In the present paper we generalize these probability distributions by introducing two new inputs which are:
\begin{enumerate}
\item The error estimate (\ref{estimation_error}) is enriched by the \emph{a priori} estimates (\ref{estimation_a_priori}) and (\ref{estimation_a_priori_discret}) to finally consider (\ref{Splitted_Error_Estimate_Approx}).
\item The probabilistic approach we develop is an extension of those we considered in \cite{ChasAs18} and \cite{arXiv_Wmp}. More precisely, by the help of an \emph{ad-hoc} density probability, we derive the probability distribution of a suitable random variable so that we can compare the two approximation errors $\|u-u^{(k_1)}_h\|_{m,p,\Omega}$ and $\|u-u^{(k_2)}_h\|_{m,p,\Omega}, (k_1 < k_2)$, considered as random variables, as will be introduced now.
\end{enumerate}
%
%
\subsection{Random solution and random approximations of determinist partial differential equation}\label{FEM_accuracy}
\noindent The purpose of this section is based of the following fundamental remark: Solution $u$ to the variational problem \textbf{(VP)}, except for particular cases, is totally unknown (being impossible to calculate it analytically); this motivates the numerical schemes one will choose to implement.\sa
This inability to determine, in most cases, the exact solution $u$, is mainly due to the complexity of the involved PDE's operator; it indeed depends on complex combinations of integrals, partial derivatives and boundary conditions, as well as on the bent geometrical shape of the domain of integration $\Omega$. All of these ingredients hence participate in the incapability to analytically determine the exact solution $u$ as their relationship with $u$ is inextricable {\sout , and then, unknown}.\sa
As a consequence, this lack of knowledge and information regarding the dependency between these ingredients and the solution $u$ motivates us to consider $u$ as a random variable, as well as any function of $u$. This paper is dedicated to the $W^{m,p}$ approximation error of $u-u^{(k)}_h$, considered as a random variable.\sa
In this frame, we view solution $u$ to the variational formulation \textbf{(VP)} defined by (\ref{VP}) in the same way as it is usual to consider the trajectory and the contact point with the ground of any solid body which is thrown, {\it i.e.}, as random. Indeed, in this case, due to the lack of information concerning the initial conditions of the trajectory of the body, the solution of the concerned inverse kinematic operator is inaccessible, and is thus seen as a random variable.\sa
In the case analyzed in this paper, the situation is much worse. Indeed, we investigate a general variational formulation \textbf{(VP)} where the analog of the inverse kinematic operator is too complex to enable us to analytically determine the corresponding solution of \textbf{(VP)} by any mathematical expression. It is one of the reasons which motivate us to view solution $u$ and its approximation $u^{(k)}_h$ as random variables, since the corresponding approximate operator conserves the complexity of the original one described above.
%
%
\subsection{The probabilistic distribution of the relative finite elements accuracy}
\noindent In the previous section we motivated the reason why we consider solution $u$ and its approximation $u^{(k)}_h$ as a random variables. \sa
To complete the description of the randomness feature of the approximation error $\|u-u^{(k)}_h\|_{m,p,\Omega}$, we also remark that, since the way a given mesh grid generator will produce any mesh is random, then the corresponding approximation $u^{(k)}_h$ is random too. \sa
For all of these reasons, based on the error estimate (\ref{Error_Estimate_Approx})-(\ref{Lambda}), we can only affirm that the value of the approximation error $\|u^{(k)}_h-u\|_{m,\Omega}$ is somewhere within the interval $[0,\beta_k]$.\sa
As a consequence, we decide to see norm $\|u-u^{(k)}_h\|_{m,p,\Omega}$ as a random variable defined as follows: \sa
Let $k, m$ and $p$ be fixed. We introduce the random variable $X_{m,p}^{(k)}$ defined by:
\begin{eqnarray}
X_{m,p}^{(k)} : & W^{m,p}(\Omega)\times W^{(k)}_{h} & \hspace{0.1cm}\longrightarrow \hspace{0.2cm}[0,\beta_k] \label{Xk_1}\noindent \label{Def_Xi_h_0} \\
& \boldsymbol{\omega}\equiv (u,u^{(k)}_h) & \hspace{0.1cm} \longmapsto \hspace{0.2cm}\D X_{m,p}^{(k)}(\boldsymbol{\omega}) = X_{m,p}^{(k)}(u,u^{(k)}_h) = \|u-u^{(k)}_h\|_{m,p,\Omega}. \label{Def_Xi_h}
\end{eqnarray}
Thus, the space product $W^{m,p}(\Omega)\times W^{(k)}_{h}$ plays the role of the usual probability space introduced in this context.\sa
Now, regarding the absence of information concerning the more likely or less likely values of norm $\|u-u^{(k)}_h\|_{m,p,\Omega}$ within the interval $[0, \beta_k]$, we will assume that the random variable $X_{m,p}^{(k)}$ is uniformly distributed over the interval $[0, \beta_k]$, with the following meaning:
\begin{equation}\label{Xk_uniform}
\forall (a,b)\in\R^2, 0 \leq a < b \leq \beta_k: Prob\left\{X_{m,p}^{(k)} \in [a,b]\right\}=\frac{b-a}{\beta_k}.
\end{equation}
Equation (\ref{Xk_uniform}) means that if one slides the interval $[a,b]$ anywhere in $[0, \beta_k]$, the probability of event $\D\left\{X_{m,p}^{(k)} \in [a,b]\right\}$ does not depend on where the interval $[a,b]$ is located in $[0, \beta_k]$, but only on its length; this corresponds to the property of uniformity of the random variable $X_{m,p}^{(k)}$. \sa
Let us now consider two families of Lagrange finite elements $P_{k_1}$ and $P_{k_2}$ corresponding to a set of values $(k_1,k_2)\in \N^2$ such that $0 < k_1 < k_2$.\sa
The two corresponding inequalities given by (\ref{Error_Estimate_Approx})-(\ref{Lambda}), assuming that solution $u$ to \textbf{(VP)} belong to $H^{k_2+1}(\Omega)$, are:
\vspace{-0.2cm}
\begin{eqnarray}
X_{m,p}^{(k_1)} \, \equiv \, \|u-u^{(k_1)}_h\|_{m,p,\Omega} & \leq & \beta_{k_1}, \label{Constante_01} \\
X_{m,p}^{(k_2)} \, \equiv \, \|u-u^{(k_2)}_h\|_{m,p,\Omega} & \leq & \beta_{k_2}, \label{Constante_02}
\end{eqnarray}
where $u^{(k_1)}_h$ and $u^{(k_2)}_h$ respectively denote the $P_{k_1}$ and $P_{k_2}$ Lagrange finite element approximations of $u$ and $\beta_{k_i}, (i=1,2),$ as defined by (\ref{Lambda}).
\begin{remark}
If one considers a given mesh for the finite element $P_{k_2}$ which contains that of $P_{k_1}$ then, for the particular class of problems where $\textbf{(VP)}$ is equivalent to a minimization formulation $\textbf{(MP)}$ (see for example \cite{ChaskaPDE}), one can show that the approximation error of the $P_{k_2}$ finite element is always smaller than that of $P_{k_1}$, and $P_{k_2}$ is more accurate than $P_{k_1}$, for all values of the mesh size $h$.
\end{remark}
Therefore, to avoid this situation, for a given value of $h$, we consider two independent meshes built by a mesh generator for $P_{k_1}$ and $P_{k_2}$. Now, usually, in order to compare the relative accuracy between these two finite elements, one asymptotically considers inequalities (\ref{Constante_01}) and (\ref{Constante_02}) to conclude that, when $h$ goes to zero, $P_{k_2}$ is more accurate that $P_{k_1}$, since $h^{k_2}$ goes faster to zero than $h^{k_1}$. \sa
However, for a given application, $h$ has a given and fixed value, so this way of comparison is not valid anymore. For this reason, our purpose is to determine the relative accuracy between two finite elements $P_{k_1}$ and $P_{k_2}, (k_1<k_2)$ for a fixed value of $h$ corresponding to two independent meshes.\sa
Moreover, since we chose to consider the two random variables $X_{m,p}^{(k_i)}, (i=1,2),$ as uniformly distributed on their respective interval of values $[0,\beta_{k_i}],(i=1,2)$, we also assume that they are independent. This assumption is, once again, the result of the lack of information which lead us to model the relationship between these two variables as independent, since any knowledge is available to more precisely localize the value of $\|u-u^{(k_1)}_h\|_{m,p,\Omega}$ if the value of $\|u-u^{(k_2)}_h\|_{m,p,\Omega}$ was known, and \emph{vice versa}.\sa
By the following result, we establish the density of probability of the random variable $Z$ defined by: $Z=X_{m,p}^{(k_1)}-X_{m,p}^{(k_2)}$.
\begin{theorem}\label{theorem_density_fZ}
Let $X_{m,p}^{(k_i)},(i=1,2),$ be the two uniform and independent random variables defined by (\ref{Def_Xi_h_0})-(\ref{Def_Xi_h}): $X_{m,p}^{(k_i)} \thicksim~U([0,\beta_{k_i}])$, where $\beta_{k_i}$ is defined by (\ref{Lambda}).\sa
Then, the random variable $Z$ defined on $\R$ has the following density of probability:
\begin{flushleft}
$f_ {Z}(z) = \hspace{0.1cm}\left| \mbox{\hfill \begin{minipage}[h]{15cm}
\begin{eqnarray}
\hspace{-0.8cm}\D 0 & \mbox{ if } & \left(z \leq -\beta_{k_2} \, \mbox{ or }\, z \geq \beta_{k_1}\right), \label{fZ_1}\\[0.2cm]
\hspace{-0.8cm}\D \frac{\beta_{k_2} + z}{\beta_{k_1}\beta_{k_2}} & \mbox{ if } & \left(\beta_{k_1} \leq \beta_{k_2}: -\beta_{k_2} \leq z \leq \beta_{k_1}-\beta_{k_2}\right) \mbox{ or } \left(\beta_{k_1} \geq \beta_{k_2}: -\beta_{k_2} \leq z \leq 0\right), \label{fZ_2}\\[0.2cm]
\hspace{-0.8cm}\D \frac{1}{\beta_{k_2}} & \mbox{ if } & \left(\beta_{k_1} \leq \beta_{k_2}: \beta_{k_1}-\beta_{k_2} \leq z \leq 0\right), \label{fZ_3}\\[0.2cm]
\hspace{-0.8cm}\D \frac{1}{\beta_{k_1}} & \mbox{ if } & \left(\beta_{k_1} \geq \beta_{k_2}: 0 \leq z \leq \beta_{k_1}-\beta_{k_2}\right), \label{fZ_4}\\[0.2cm]
\hspace{-0.8cm}\D \frac{\beta_{k_1} - z}{\beta_{k_1}\beta_{k_2}} & \mbox{ if } & \left(\beta_{k_1} \leq \beta_{k_2}: 0 \leq z \leq \beta_{k_1}\right) \mbox{ or } \left(\beta_{k_1} \geq \beta_{k_2}: \beta_{k_1}-\beta_{k_2} \leq z \leq \beta_{k_1}\right). \label{fZ_5} 
\end{eqnarray}
\end{minipage}
}
\right.$
\end{flushleft}
\end{theorem}
\begin{prooff}
Let us now remark that, since the support of the two random variables $X_{m,p}^{(k_i)}, (i=1,2),$ is $[0,\beta_{k_i}]$, the support of the density $f_Z$ is therefore $[-\beta_{k_2},\beta_{k_1}]$, which corresponds to (\ref{fZ_1}).\sa
Let us consider the case where $\beta_{k_1} \leq \beta_{k_2}$. \sa
If $f_{X_{m,p}^{(k_i)}}(x_i)$ denotes the density of probability defined on $\R$ associated to each random variable $X_{m,p}^{(k_i)}$, since we assume they are independent variables, the density $f_Z(z)$ is given by:
\begin{equation}\label{densite_fZ_integrale}
\D f_Z(z) = \int_{\R}f_{X_{m,p}^{(k_2)}}(x_2)f_{X_{m,p}^{(k_1)}}(x_2+z)dx_2,
\end{equation}
where:
\begin{equation}\label{densite_fX}
\D f_{X_{m,p}^{(k_i)}}(x_i) = \frac{1}{\beta_{k_i}}\mathbbm{1}_{[0,\beta_{k_i}]}(x_i), (i=1,2),
\end{equation}
and $\mathbbm{1}_{[c,d]}$ is the indicator function of the interval $[c,d], \forall (c,d)\in \R^2$.\sa
Furthermore, due to the definition (\ref{densite_fX}) of the density for each variable $X_{m,p}^{(k_i)}, (i=1,2),$ the integrand of (\ref{densite_fZ_integrale}) can be expressed as follows:
\begin{flushleft}
$\hspace{0.75cm} f_{X_{m,p}^{(k_2)}}(x_2)f_{X_{m,p}^{(k_1)}}(x_2+z) = \hspace{0.1cm}\left| \mbox{\hfill \begin{minipage}[h]{11.5cm}
\begin{eqnarray}
\hspace{-2.2cm}\D \frac{1}{\beta_{k_1}\beta_{k_2}} & \mbox{ if } & 0 \leq x_2 \leq \beta_{k_2} \,  \mbox{ and } \, 0 \leq x_2+z \leq \beta_{k_1}, \medskip \\
\hspace{-2.2cm}\D 0 & \mbox{if } & \mbox{not}.
\end{eqnarray}
\end{minipage}
}
\right.$
\end{flushleft}
\vspace{0.15cm}
As a consequence, $x_2 \in [0,\beta_{k_2}] \cap [-z,-z+\beta_{k_1}]$, which leads one to consider the five following cases corresponding to the significant relative positions between the intervals $[0,\beta_{k_2}]$ and $[-z,-z+\beta_{k_1}]$ .
\begin{enumerate}
\item Let us assume that $-z+\beta_{k_1} \leq 0$. If $z\geq \beta_{k_1}$ then $[0,\beta_{k_2}] \cap [-z,-z+\beta_{k_1}] = \varnothing$ and $f_Z(z)=0$ which is again the result of (\ref{fZ_1}).
\item We consider now the values of $z$ such that $-z\leq 0$ and $0 \leq -z+\beta_{k_1} \leq \beta_{k_2}$. \\  If $0 \leq z \leq \beta_{k_1}$, then $x_2 \in [0,-z+\beta_{k_1}]$ and by (\ref{densite_fZ_integrale}) we get the expected expression of (\ref{fZ_5}) since we have:
\begin{equation}
\D f_Z(z) = \int_{0}^{-z+\beta_{k_1}}\frac{dx_2}{\beta_{k_1}\beta_{k_2}} = \frac{\beta_{k_1}-z}{\beta_{k_1}\beta_{k_2}}, (0 \leq z \leq \beta_{k_1}).
\end{equation}
\item The following case concerns the values of $z$ such that $-z \geq 0$ and $-z+\beta_{k_1} \leq \beta_{k_2}$. If $\beta_{k_1}-\beta_{k_2} \leq z \leq 0$, then $x_2 \in [-z,-z+\beta_{k_1}]$ and we get (\ref{fZ_3}) since:
\begin{equation}
\D f_Z(z) = \int_{-z}^{-z+\beta_{k_1}}\frac{dx_2}{\beta_{k_1}\beta_{k_2}} = \frac{1}{\beta_{k_2}}, (\beta_{k_1}-\beta_{k_2} \leq z \leq 0).
\end{equation}
\item If $-z \leq 0$ and $-z+\beta_{k_1} \geq \beta_{k_2}$ then $0\leq z \leq \beta_{k_1}-\beta_{k_2}$ which is impossible since $\beta_{k_1}\leq\beta_{k_2}$.
\item We consider now the values of $z$ such that $0 \leq -z \leq \beta_{k_2}$ and $-z+\beta_{k_1} \geq \beta_{k_2}$. \sa If $-\beta_{k_2} \leq z \leq
\beta_{k_1}-\beta_{k_2}$, then $x_2 \in [-z,\beta_{k_2}]$ and, since we have:
\begin{equation}
\D f_Z(z) = \int_{-z}^{\beta_{k_2}}\frac{dx_2}{\beta_{k_1}\beta_{k_2}} = \frac{\beta_{k_2}+z}{\beta_{k_1}\beta_{k_2}}, (-\beta_{k_2}\leq z \leq \beta_{k_1}-\beta_{k_2}),
\end{equation}
we get the expected expression of (\ref{fZ_2}).
\item Finally, if $-z \geq \beta_{k_2}$, then $[0,\beta_{k_2}] \cap [-z,-z+\beta_{k_1}] = \varnothing$ and $f_Z(z)=0$ which corresponds to (\ref{fZ_1}).
\end{enumerate}
\vs
The other cases corresponding to $\beta_{k_1} \geq \beta_{k_2}$ can be deduced using the same arguments.
\end{prooff}
From Theorem \ref{theorem_density_fZ}, one can infer the entire cumulative distribution function $F_Z(z)$ defined by:
\begin{equation}\label{fonction repartition_FZ}
F_Z(z) = \int_{-\infty}^{z}f_Z(z) dz.
\end{equation}
However, since we are interested in determining the more likely finite element between $P_{k_1}$ and $P_{k_2}, (k_1 <~k_2)$, we focus the following corrolary to the value of $F_Z(0)$ which corresponds to $Prob\left\{ X_{m,p}^{(k_1)} \leq X_{m,p}^{(k_2)}\right\}$.
\begin{corollary}
Let $X_{m,p}^{(k_i)},(i=1,2),$ be the two uniform and independent random variables defined by (\ref{Def_Xi_h_0})-(\ref{Def_Xi_h}): $X_{m,p}^{(k_i)} \thicksim~U([0,\beta_{k_i}])$, where $\beta_{k_i}$ is defined by (\ref{Lambda}). Then, we have:
\begin{flushleft}
$\hspace{1.34cm} Prob\left\{ X_{m,p}^{(k_1)} \leq X_{m,p}^{(k_2)}\right\} = \hspace{0.1cm}\left| \mbox{\hfill \begin{minipage}[h]{11cm}
\begin{eqnarray}
\hspace{-6cm}\D 1-\frac{1}{2}\frac{\beta_{k_1}}{\beta_{k_2}} & \mbox{ if } & \beta_{k_1} \leq \beta_{k_2}, \label{Prob_beta_k1_beta_k2_1} \medskip \\
\hspace{-6cm}\D \frac{1}{2}\frac{\beta_{k_2}}{\beta_{k_1}} & \mbox{if } & \beta_{k_1} \geq \beta_{k_2}. \label{Prob_beta_k1_beta_k2_2}
\end{eqnarray}
\end{minipage}
}
\right.$
\end{flushleft}
\end{corollary}
\begin{prooff}
\begin{itemize}
\item Let us consider the case where $\beta_{k_1} \leq \beta_{k_2}$. Then, by definition (\ref{fonction repartition_FZ}) of the entire cumulative distribution function $F_Z(z)$, at $z=0$, we have:
\begin{equation}\label{fonction repartition_FZ_en_z=0_1}
\D F_Z(0) = \int_{-\infty}^{0}f_Z(z) dz = \int_{-\beta_{k_2}}^{\beta_{k_1}\!-\beta_{k_2}}\frac{\beta_{k_2} + z}{\beta_{k_1}\beta_{k_2}}\,\, dz \, + \, \int_{\beta_{k_1}\!-\beta_{k_2}}^{0}\frac{dz}{\beta_{k_2}} \, = \, 1-\frac{1}{2}\frac{\beta_{k_1}}{\beta_{k_2}},
\end{equation}
where we used the values of the density $f_Z(z)$ given by (\ref{fZ_2}) and (\ref{fZ_3}).
\item In the same way, when $\beta_{k_1} \geq \beta_{k_2}$ we have:
\begin{equation}\label{fonction repartition_FZ_en_z=0_2}
\D F_Z(0) = \int_{-\infty}^{0}f_Z(z)\,dz \, = \,\int_{-\beta_{k_2}}^{0}\frac{\beta_{k_2} + z}{\beta_{k_1}\beta_{k_2}}\,\, dz \, = \, \frac{1}{2}\frac{\beta_{k_2}}{\beta_{k_1}}.
\end{equation}
\end{itemize}
\end{prooff}
We can now explicit the probability distribution of the event $\D\left\{ X_{m,p}^{(k_1)} \leq X_{m,p}^{(k_2)}\right\}$ given by (\ref{Prob_beta_k1_beta_k2_1})-(\ref{Prob_beta_k1_beta_k2_2}) as a function of the mesh size $h$. \sa
To this end, we remark that, since each $\beta_{k_i} (i=1,2)$ has two possible values depending on the relative position between $h$ and $\hslash_{k_i}$ defined by (\ref{hslash_k}), the probability distribution we are looking for must be splitted in the corresponding cases as well.
 \sa
Let us hence introduce the constants $C_{k_i}, (i=1,2),$ defined by:
\begin{equation}\label{Ck_i}
C_{k_i} = \mathscr{C}_{k_i}|u|_{k_i+1,p,\Omega},
\end{equation}
and the specific value $h_{k_1,k_2}^*$ of $h$ which corresponds to the intersection of the curves $\varphi_{k_i}, (i=~1,2),$ defined by $\varphi_{k_i}(h)\equiv C_{k_i}h^{k_i+1-m}$, (see Figure \ref{Courbes_1}).\sa
Then, we have:
\begin{equation}\label{h*k1k2}
\D h_{k_1,k_2}^* = \left(\frac{C_{k_1}}{C_{k_2}}\right)^{\frac{1}{k_2-k_1}} \equiv  \left(\frac{\mathscr{C}_{k_1}}{\mathscr{C}_{k_2}}\frac{|u|_{k_1+1,p,\Omega}}{|u|_{k_2+1,p,\Omega}}\right)^{\frac{1}{k_2-k_1}}.
\end{equation}
We notice that $h_{k_1,k_2}^*$ and $\hslash_{k_i}, (i=1,2),$ strongly depend on $m$ and $p$, since $\mathscr{C}_{k_i}$ and $\beta_{k_i}$ depend on these two parameters as well. As a consequence, in the following theorem the different formulas of $Prob\left\{ X_{m,p}^{(k_1)} \leq X_{m,p}^{(k_2)}\right\}$ will contain this dependency on $m$ and $p$.
\begin{theorem}\label{Corollaire_Loi P_k1_k2_V2}
Let $X_{m,p}^{(k_i)},(i=1,2),$ be the two uniform and independent random variables defined by (\ref{Def_Xi_h_0})-(\ref{Def_Xi_h}). Then, $\D {\cal P}_{k_1,k_2}(h) \equiv Prob\left\{ X_{m,p}^{(k_1)} \leq X_{m,p}^{(k_2)}\right\}$ is determined by:\sa
If $h_{k_1,k_2}^* \geq \max(\hslash_{k_1},\hslash_{k_2})$ then $\hslash_{k_1}\leq\hslash_{k_2}$ and:
\begin{flushleft}
\hspace{1cm}$\D {\cal P}_{k_1,k_2}(h) \, = \,\left| \mbox{\hfill \begin{minipage}[h]{12cm}
\begin{eqnarray}
\D \hspace{-5cm}\frac{1}{2}\left(\frac{h}{h_{k_1,k_2}^*}\right)^{k_2-k_1} & \mbox{ if } & 0 \leq h \leq \hslash_{k_1}, \label{F1}\\[0.2cm]
\D \hspace{-5cm}\frac{1}{2}\left(\frac{h}{\hslash_{k_2}}\right)^{k_2+1-m} & \mbox{ if } & \hslash_{k_1}\leq h \leq\hslash_{k_2}, \label{F2}\\[0.2cm]
\D \hspace{-5cm}\frac{1}{2}\hspace{2.45cm}& \mbox{ if } & h \geq \hslash_{k_2}. \label{F3}
\end{eqnarray}
\end{minipage}
}
\right.$
\end{flushleft}
If $h_{k_1,k_2}^* \leq \min(\hslash_{k_1},\hslash_{k_2})$ then $\hslash_{k_2}\leq\hslash_{k_1}$ and:
\begin{flushleft}
\hspace{1cm}$\D {\cal P}_{k_1,k_2}(h) \, = \,\left| \mbox{\hfill \begin{minipage}[h]{12cm}
\begin{eqnarray}
\D \hspace{-3.5cm}\frac{1}{2}\left(\frac{h}{h_{k_1,k_2}^*}\right)^{k_2-k_1} & \mbox{ if } & 0 \leq h \leq h_{k_1,k_2}^*, \label{F4}\\[0.2cm]
\D \hspace{-3.5cm}1-\frac{1}{2}\left(\frac{h_{k_1,k_2}^*}{h}\right)^{k_2-k_1} & \mbox{ if } & h_{k_1,k_2}^* \leq h \leq \hslash_{k_2}, \label{F5}\\[0.2cm]
\D \hspace{-3.5cm}1-\frac{1}{2}\left(\frac{h}{\hslash_{k_1}}\right)^{k_1+1-m} & \mbox{ if } & \hslash_{k_2}\leq h \leq\hslash_{k_1}, \label{F6}\\[0.2cm]
\D \hspace{-3.5cm}\frac{1}{2}\hspace{2.45cm}& \mbox{ if } & h \geq \hslash_{k_1}. \label{F7}
\end{eqnarray}
\end{minipage}
}
\right.$
\end{flushleft}
\end{theorem}
\begin{prooff}
To establish the proof of Theorem \ref{Corollaire_Loi P_k1_k2_V2}, we will consider a geometrical interpretation of the error estimate (\ref{Constante_01}) and (\ref{Constante_02}). \sa
These two inequalities can indeed be geometrically viewed with the help of the relative position between the two curves $\varphi_{k_i}(h), (i=1,2),$ introduced above, and the horizontal line defined by $\D \psi(h)\equiv\frac{\|l\|_{V'}}{\alpha^*}$, (see Figure \ref{Courbes_1}).\sa
Then, depending on the position of the horizontal line $\psi(h)$ with the particular value $C_{k_i}h_{k_1,k_2}^*, (i=1,2),$ we have to consider the two following cases:
\begin{itemize}
\item If $\D\frac{\|l\|_{V'}}{\alpha^*}\leq C_{k_1}h_{k_1,k_2}^* = C_{k_2}h_{k_1,k_2}^*$, then $\hslash_{k_1}\leq\hslash_{k_2}\leq h_{k_1,k_2}^*$.
\item If $\D\frac{\|l\|_{V'}}{\alpha^*}\geq C_{k_1}h_{k_1,k_2}^* = C_{k_2}h_{k_1,k_2}^*$, then $h_{k_1,k_2}^*\leq\hslash_{k_2}\leq\hslash_{k_1}$.
\end{itemize}
\begin{figure}[h]
  \centering
  \includegraphics[width=14cm]{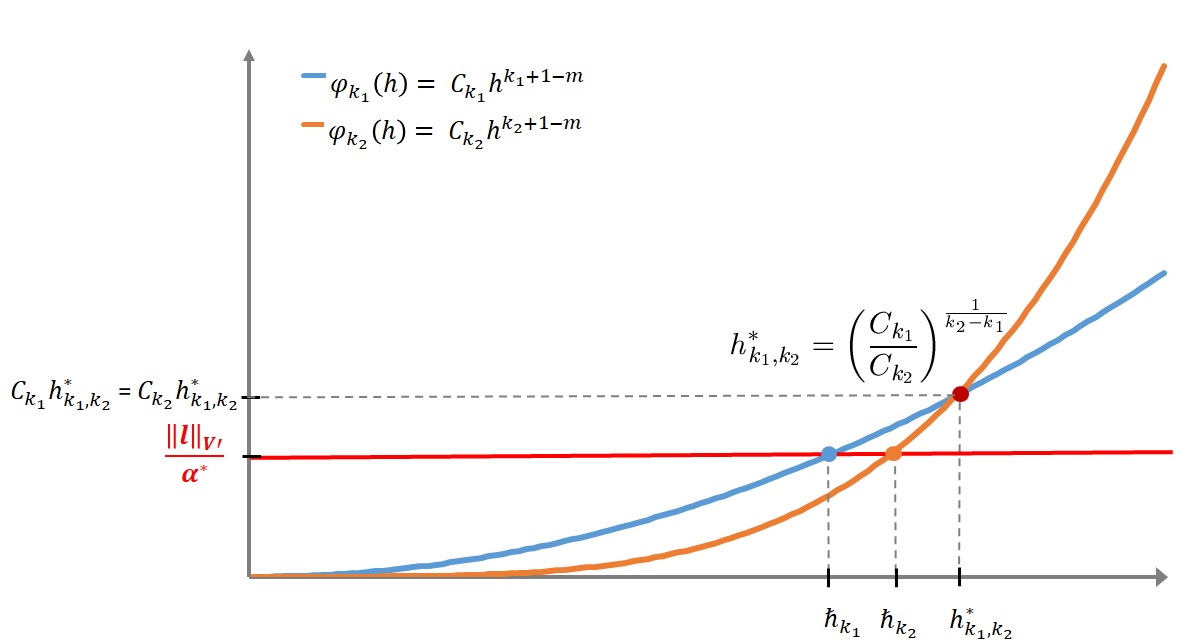}
  \caption{Relative positions between the curves $C_{k_i}h^{k_i+1-m}$ and the line $\|l\|_{V'} / \alpha^*$.} \label{Courbes_1}
\end{figure}
So, let us consider the first case when $\D\frac{\|l\|_{V'}}{\alpha^*}\leq C_{k_1}h_{k_1,k_2}^* = C_{k_2}h_{k_1,k_2}^*$, or equivalently, $h_{k_1,k_2}^* \geq \max(\hslash_{k_1},\hslash_{k_2})$. Therefore, due to the relative positions between the three curves $\varphi_{k_1}, \varphi_{k_2}$ and $\psi$, we have the following results:
\begin{enumerate}
\item $\D\forall h \leq \hslash_{k_1}: C_{k_2}h^{k_2+1-m} \leq C_{k_1}h^{k_1+1-m} \leq \frac{\|l\|_{V'}}{\alpha^*}$ and $\beta_{k_2}\leq\beta_{k_1}$  with $\beta_{k_i}=C_{k_i}h^{k_{i}+1-m}, (i=~1,2)$. Then, from (\ref{Prob_beta_k1_beta_k2_2}), we get (\ref{F1}).\vspace{0.1cm}
\item $\D\forall \hslash_{k_1} \leq h \leq \hslash_{k_2}: C_{k_2}h^{k_2+1-m} \leq \frac{\|l\|_{V'}}{\alpha^*} \leq C_{k_1}h^{k_1+1-m}$ and $\beta_{k_2}\leq\beta_{k_1}$ with $\D\beta_{k_1}=\frac{\|l\|_{V'}}{\alpha^*}$ and $\beta_{k_2}=C_{k_2}h^{k_{2}+1-m}$. Then, again from (\ref{Prob_beta_k1_beta_k2_2}), we get (\ref{F2}).\vspace{0.1cm}
\item $\D\forall \hslash_{k_2} \leq h \leq h_{k_1,k_2}^*: \frac{\|l\|_{V'}}{\alpha^*} \leq C_{k_2}h^{k_2+1-m} \leq  C_{k_1}h^{k_1+1-m}$ and $\D\beta_{k_1}=\beta_{k_2}=\frac{\|l\|_{V'}}{\alpha^*}$. Then, from (\ref{Prob_beta_k1_beta_k2_1}) or (\ref{Prob_beta_k1_beta_k2_2}), we get (\ref{F3}).\vspace{0.1cm}
\item $\D\forall h_{k_1,k_2}^* \leq h : \frac{\|l\|_{V'}}{\alpha^*} \leq C_{k_1}h^{k_1+1-m} \leq  C_{k_2}h^{k_2+1-m}$ and $\D\beta_{k_1}=\beta_{k_2}=\frac{\|l\|_{V'}}{\alpha^*}$. Then, again from (\ref{Prob_beta_k1_beta_k2_1}) or (\ref{Prob_beta_k1_beta_k2_2}), we also get (\ref{F3}).
\end{enumerate}
We now consider the case where $\D\frac{\|l\|_{V'}}{\alpha^*}\geq C_{k_1}h_{k_1,k_2}^* = C_{k_2}h_{k_1,k_2}^*$ which corresponds to $h_{k_1,k_2}^* \leq \min(\hslash_{k_1},\hslash_{k_2})$. Therefore, using once again the relative positions between the curves $\varphi_{k_1}, \varphi_{k_2}$ and $\psi$, we deduce the following results:
\begin{enumerate}
\item $\D\forall h \leq h_{k_1,k_2}^*: C_{k_2}h^{k_2+1-m} \leq C_{k_1}h^{k_1+1-m} \leq \frac{\|l\|_{V'}}{\alpha^*}$ and $\beta_{k_2}\leq\beta_{k_1}$  with $\beta_{k_i}=C_{k_i}h^{k_{i}+1-m}, (i=1,2)$. Then, from (\ref{Prob_beta_k1_beta_k2_2}), we get (\ref{F4}).\vspace{0.1cm}
\item $\D\forall h_{k_1,k_2}^* \leq h \leq \hslash_{k_2}: C_{k_1}h^{k_1+1-m} \leq C_{k_2}h^{k_2+1-m} \leq \frac{\|l\|_{V'}}{\alpha^*}$ and $\beta_{k_1}\leq\beta_{k_2}$ with $\beta_{k_i}=C_{k_i}h^{k_{i}+1-m}, \\(i=1,2)$. Then, from (\ref{Prob_beta_k1_beta_k2_1}), we get (\ref{F5}).\vspace{0.1cm}
\item $\D\forall \hslash_{k_2} \leq h \leq \hslash_{k_1}: C_{k_1}h^{k_1+1-m} \leq \frac{\|l\|_{V'}}{\alpha^*} \leq C_{k_2}h^{k_2+1-m}$ and $\beta_{k_1}\leq\beta_{k_2}$ with $\beta_{k_1}=C_{k_1}h^{k_1+1-m}$ and $\D\beta_{k_2}=\frac{\|l\|_{V'}}{\alpha^*}$. Then, again from (\ref{Prob_beta_k1_beta_k2_1}), we get (\ref{F6}).\vspace{0.1cm}
\item $\D\forall h \geq \hslash_{k_1} : \frac{\|l\|_{V'}}{\alpha^*} \leq C_{k_1}h^{k_1+1-m} \leq  C_{k_2}h^{k_2+1-m}$ and $\D\beta_{k_1}=\beta_{k_2}=\frac{\|l\|_{V'}}{\alpha^*}$. Then, from (\ref{Prob_beta_k1_beta_k2_1}) or (\ref{Prob_beta_k1_beta_k2_2}), we get (\ref{F7}).
\end{enumerate}
\end{prooff}
\begin{remark}
We notice that the probability distribution ${\cal P}_{k_1,k_2}(h)$ given by Theorem \ref{Corollaire_Loi P_k1_k2_V2} generalizes those we found in \cite{ChasAs18} and \cite{arXiv_Wmp}. Indeed, when we derived the probability distribution ${\cal P}_{k_1,k_2}(h)$ without taking into account the a priori estimates (\ref{estimation_a_priori}) and (\ref{estimation_a_priori_discret}), we obtained:
\begin{flushleft}
$\D \hspace{1.5cm}{\cal P}_{k_1,k_2}(h) \, = \,\left| \mbox{\hfill \begin{minipage}[h]{8cm}
\begin{eqnarray}
\hspace{-0.5cm}\D \frac{1}{2}\!\left(\!\frac{\!\!h}{h_{k_1,k_2}^*}\!\right)^{\!\!k_2-k_1} & \hspace{-0.2cm}\mbox{ if }\hspace{-0.2cm} & \hs 0 \leq h \leq h_{k_1,k_2}^*, \label{CMAM1}\\[0.2cm]
\hspace{-0.5cm}\D 1 - \frac{1}{2}\!\left(\!\frac{h_{k_1,k_2}^*}{\!\!h}\!\right)^{\!\!k_2-k_1} & \hspace{-0.2cm}\mbox{ if }\hspace{-0.2cm} & h \geq h_{k_1,k_2}^*.\label{CMAM2}
\end{eqnarray}
\end{minipage}
}
\right.$
\end{flushleft}
However, in Theorem \ref{Corollaire_Loi P_k1_k2_V2} the contribution of the a priori estimates (\ref{estimation_a_priori}) and (\ref{estimation_a_priori_discret}) modifies the probability distribution given by (\ref{CMAM1})-(\ref{CMAM2}), since the horizontal line $\D \psi(h)=\frac{\|l\|_{V'}}{\alpha^*}$ interferes with the two polynomials $C_{k_i}h^{k_i+1-m}, (i=1,2)$, (see Figure \ref{Courbes_1}).
\end{remark}
This phenomenon will be analyzed in the following section.
\section{Discussion}\label{discussion}
\noindent First of all, let us plot the shape of the two cases of the probability distribution we derived in Theorem \ref{Corollaire_Loi P_k1_k2_V2} and those we got in \cite{ChasAs18} and \cite{arXiv_Wmp} that we recalled in (\ref{CMAM1})-(\ref{CMAM2}).
\begin{figure}[h]
  \centering
  \includegraphics[width=14cm]{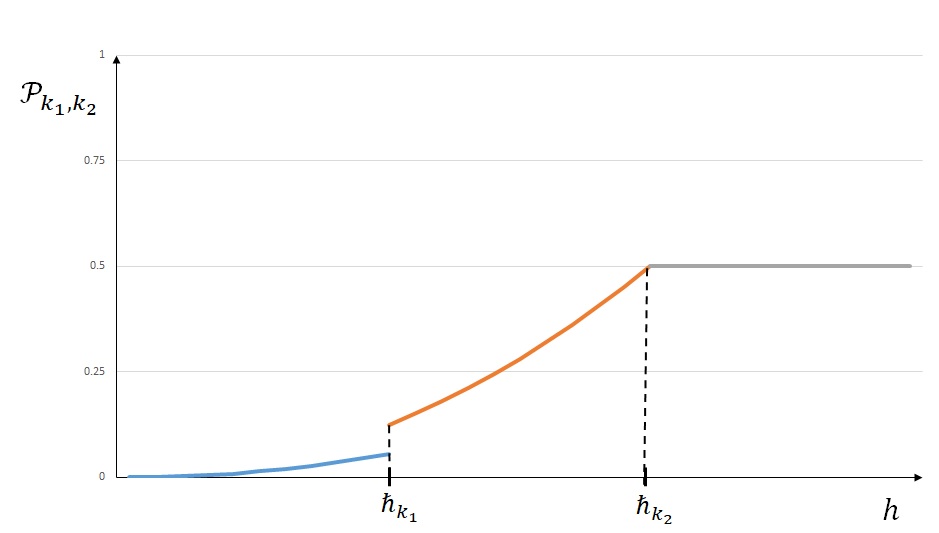}
  \caption{Probability distribution ${\cal P}_{k_1,k_2}(h)$ when $h_{k_1,k_2}^* \geq \max(\hslash_{k_1},\hslash_{k_2})$.} \label{Corrolaire_Case_1}
\end{figure}
\begin{figure}[h]
  \centering
  \includegraphics[width=14cm]{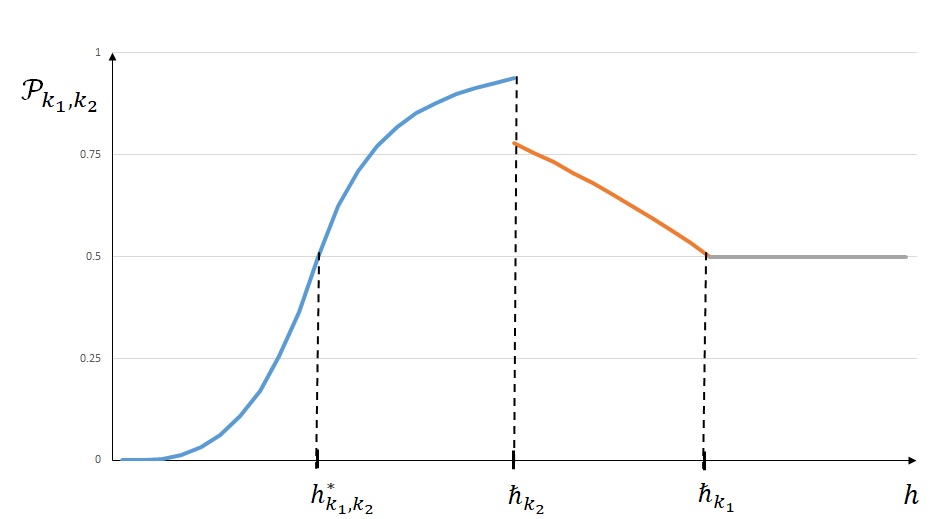}
  \caption{Probability distribution ${\cal P}_{k_1,k_2}(h)$ when $h_{k_1,k_2}^* \leq \min(\hslash_{k_1},\hslash_{k_2})$.} \label{Corrolaire_Case_2}
\end{figure}
\begin{figure}[h]
  \centering
  \includegraphics[width=14cm]{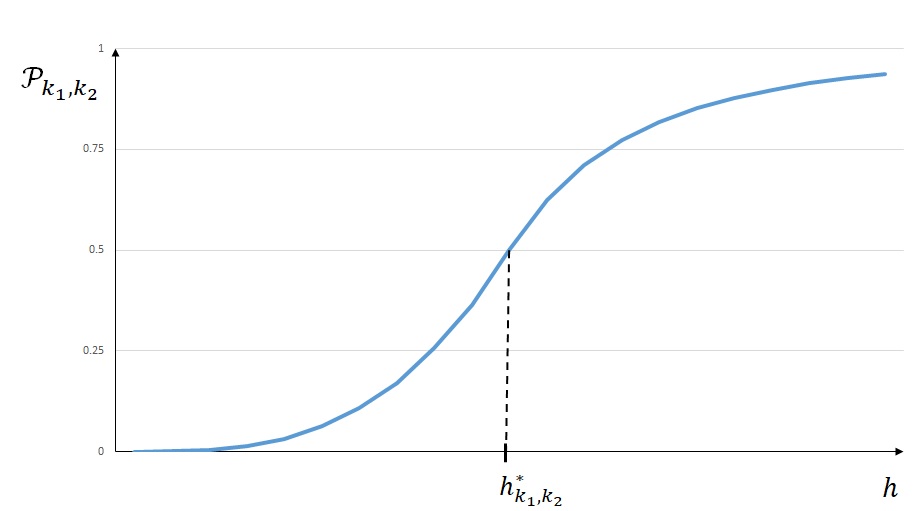}
  \caption{Probability distribution ${\cal P}_{k_1,k_2}(h)$ published in \cite{ChasAs18}.} \label{CMAM}
\end{figure}
As one can see, the distribution (\ref{CMAM1})-(\ref{CMAM2}) plotted in Figure \ref{CMAM} resembles more the second case of Theorem \ref{Corollaire_Loi P_k1_k2_V2} plotted in Figure \ref{Corrolaire_Case_2} than its first case, as plotted in Figure \ref{Corrolaire_Case_1}. \sa
In fact, in the latter case, when $\hslash_{k_1}\leq\hslash_{k_2}$, the line $\D \psi(h)=\frac{\|l\|_{V'}}{\alpha^*}$ interacts with the two polynomials $C_{k_i}h^{k_i+1-m}, (i=1,2),$ before the critical value $h^*_{k_1,k_2}$. As a consequence, its contribution correspondingly takes place in the expression of the probability distribution for $h\geq\hslash_{k_1}$. \sa
On the contrary, when $\hslash_{k_1}\geq\hslash_{k_2}$, the contribution of the line $\D \psi(h)=\frac{\|l\|_{V'}}{\alpha^*}$ has only to be taken into account after the value of $h^*_{k_1,k_2}$. But, in this case, when $h\geq\hslash_{k_2}$, once again line $\psi(h)$ interacts with the two polynomials $C_{k_i}h^{k_i+1-m}, (i=1,2),$ and the equivalent distribution we found in (\ref{CMAM1})-(\ref{CMAM2}) is cut to give the corresponding one in (\ref{F6})-(\ref{F7}) plotted in Figure \ref{Corrolaire_Case_2}.\sa
This shows the importance of the \emph{a priori} estimates (\ref{estimation_a_priori}) and (\ref{estimation_a_priori_discret}) which are not considered in the classical point of view limited to the asymptotic behavior of the error estimate (\ref{estimation_error}), when the mesh size $h$ goes to zero. The reason for this is the fact that the right hand sides of (\ref{estimation_a_priori}) and (\ref{estimation_a_priori_discret}) do not depend on $h$; it does hence bring no information for the desired asymptotic behavior.\sa
The goal of the probabilistic approach we developed in this paper is to evaluate, for any \underline{fixed} value of the mesh size $h$, the relative accuracy between two Lagrange finite elements $P_{k_1}$ and $P_{k_2}, (k_1\leq k_2)$. The probability distributions of Theorem \ref{Corollaire_Loi P_k1_k2_V2} adjust those we had found in \cite{ChasAs18} and \cite{arXiv_Wmp}. \sa
Indeed, probability distribution (\ref{CMAM1})-(\ref{CMAM2}) claims that finite element $P_{k_1}$ is more likely accurate than the $P_{k_2}$ one with a high level of probability when $h\geq h^*_{k_1,k_2}$. Actually, this probability can reach the value 1, as it can be observed in Figure \ref{CMAM} and $P_{k_1}$ could be \emph{almost surely} more accurate than $P_{k_2}$. So, distributions (\ref{F1})-(\ref{F3}) and (\ref{F4})-(\ref{F7}) limit the value of this probability.\sa
More precisely, if $\hslash_{k_1}\leq\hslash_{k_2}$, then the probability such that $P_{k_1}$ is more likely accurate than $P_{k_2}$ will never be greater than 0.5 (see Figure \ref{Corrolaire_Case_1}), and if $\hslash_{k_1}\geq\hslash_{k_2}$, then for $h\geq\hslash_{k_2}$ the corresponding probability will decrease to the limit value of 0.5, even if it increases with values greater than 0.5 for $h$ values between $h^*_{k_1,k_2}$ and $\hslash_{k_2}$, (see Figure \ref{Corrolaire_Case_2}).
\section{Conclusions}\label{Conclusion}
\subsection{The relative accuracy between two finite elements}
\noindent In this paper, we presented a generalized probability distribution which takes into account the three available standard estimates one can deal with a solution $u$ to a variational formulation together with its approximation $u_{h}^{(k)}$ performed by Lagrange finite elements $P_k$. \sa
Those are the \emph{a priori} estimates for the solution $u$ and its approximation $u_{h}^{(k)}$ carried out by finite elements, on the one hand, and the error estimate, on the other hand. \sa
More precisely, by taking into account the \emph{a priori} estimates (\ref{estimation_a_priori}) and (\ref{estimation_a_priori_discret}), we enriched the probability distribution (\ref{CMAM1})-(\ref{CMAM2}) we got in \cite{ChasAs18} and \cite{arXiv_Wmp} to find the new one (\ref{F1})-(\ref{F7}) which was derived in Theorem \ref{Corollaire_Loi P_k1_k2_V2}.\sa
As we saw in the previous paragraph, if the contribution of the \emph{a priori} estimates (\ref{estimation_a_priori}) and (\ref{estimation_a_priori_discret}) are not considered when one limits the analysis of the relative accuracy between two finite elements to the asymptotic rate of convergence when the mesh size $h$ goes to zero, one must deal with these two estimates to conclude about the more likely accurate one when $h$ has a fixed value.\sa
Moreover, regarding the probabilistic approach we developed in this paper, the \emph{a priori} estimates (\ref{estimation_a_priori}) and (\ref{estimation_a_priori_discret}) brought a more realistic behavior of the distribution (\ref{F1})-(\ref{F7}) in comparison with (\ref{CMAM1})-(\ref{CMAM2}). Indeed, the probability given by (\ref{F1})-(\ref{F7}) cannot anymore reach a value of 1 since it was quite a surprise in (\ref{CMAM1})-(\ref{CMAM2}) to obtain this theoretical value. This  result meant that the event ``$P_{k_1}$ \emph{is more accurate than} $P_{k_2}$" is an \emph{almost} sure event when $k_1<k_2$.
\subsection{Generalization to other numerical methods}
\noindent Finally, we would like to mention that the probabilistic approach we proposed in this work is not restricted to the finite elements method but may be extended to other types of approximation: given a class of numerical schemes and their corresponding approximation error estimates, one is able to order them, not only in terms of asymptotic rate of convergence, but also by evaluating the most probably accurate. \sa
It is for example the case of the numerical integration where the composite quadrature error has a mathematical structure which looks like the error estimate (\ref{estimation_error}) we considered in the present work. \sa
More precisely, as an example, for a composite quadrature of order $k$ on an a given interval $[a,b]$, if $f\in C^{k+1}([a,b])$ is a given function, the corresponding composite quadrature error can be written \cite{Crouzeix} as:
\begin{equation} \label{erreur_ quadrature}
\D \left|\int_{a}^{b}f(x)dx - \sum_{i=0}^{N}\lambda_i f(x_i)\right| \, \leq \, C_k\,h^{k+1},
\end{equation}
where $h$ denotes the size of the equally spaced panels which discretized the interval $[a,b]$, $(\lambda_i,x_i), (i=0,N),$ are $(2N+2)$ given numbers such that all the $x_i$ belong to $[a,b]$, and $C_k$ is a constant independent of $h$ which mainly depends on $f$ and $k$.\sa
As a consequence, due to the similar mathematical structure between (\ref{erreur_ quadrature}) and (\ref{estimation_error}), with the same arguments we introduced to compare the relative accuracy between two Lagrange finite elements $P_{k_1}$ and $P_{k_2}, (k_1<k_2),$ one can evaluate the probability of the more accurate numerical composite quadrature associated to two different parameters $k_1$ and $k_2, (k_1\leq k_2)$, for a \underline{fixed} value of $h$.\sa
This will make sense because, usually, to bypass the lack of information associated with the unknown value of the left hand side of (\ref{erreur_ quadrature}) in the interval $[0, C_k h^{k+1}]$, only the asymptotic convergence rate comparison is concerned to appreciate the relative accuracy between the two numerical quadratures of order $k_1$ and $k_2, (k_1<k_2)$.\sa
Nevertheless, this procedure is no longer available when one wants to compare two composite quadratures in the case when the size of the equally panels $h$ is fixed, as it is for any application. Thus, the probabilistic approach we propose here could be a relevant alternative. \sa
The same consideration may be developed to compare the accuracy between two numerical schemes of order $k_1$ and $k_2$ which are performed to approximate the exact solution of an ordinary differential equations. To precise these ideas, let us consider the solution $y$ of the first order ordinary differential initial value problem defined on a given interval $[t_0,t_0+T]$:
\begin{equation}\label{CP}
\textbf{(CP}\textbf{)} \hspace{0.2cm} \left\{
\begin{array}{lll}
y'(t) & = & f(t,y(t)), \, t\in[t_0,t_0+T], \\[0.2cm]
y(t_0) & = & y_0,
\end{array}
\right.
\end{equation}
where $y_0$ is given. \sa
Let us also restrict ourselves by considering one-step numerical methods to approximate the function $y$, solution to problem \textbf{(CP)}. \sa
Namely, if we introduce a constant mesh size $h=t_{n+1}-t_{n}$, where $(t_n)_{n=0,N}$ denotes the sequence of $(N+1)$ values of $t$ within the interval $[t_0, t_0+T]$, then the corresponding one-step numerical scheme is given by:
\begin{equation}\label{CP_Approx}
\textbf{(CP)$_h$} \hspace{0.2cm} \left\{
\begin{array}{lll}
y_{n+1} & = & y_{n} + h\Phi(t_n,y_n), \, n>0, \\[0.2cm]
y_0 & = & y_0,
\end{array}
\right.
\end{equation}
where $\Phi$ is a given function "sufficiently" smooth which characterizes the numerical scheme \textbf{(CP)}$_h$. \sa
Moreover, \textbf{(CP)}$_h$ is called a numerical scheme of order $k$ if \cite{Crouzeix}:
\begin{equation}\label{erreur_EDO}
\D \sum_{n=0}^{N-1}|y(t_{n+1})-y(t_{n})- h\Phi(t_n,y(t_{n}))| \leq C_k\,h^k,
\end{equation}
where $C_k$ is a constant independent of $h$ which depends on $y$, $\Phi$ and $k$, (see \cite{Stroud} for the dependency on $k$). \sa
So, when considering two one-step numerical methods of order $k_1$ and $k_2, (k_1<k_2),$ defined by two functions $\Phi_{1}$ and $\Phi_{2}$, and due to the similar structure between (\ref{erreur_EDO}) and (\ref{estimation_error}), one would be able to evaluate the probability of the more accurate scheme with the same arguments we implemented when comparing the relative accuracy between $P_{k_1}$  and $P_{k_2}$ Lagrange finite elements. \sa
In summary, when one wants to evaluate the relative accuracy between two numerical methods which belong to a given family of approximations, the probabilistic approach we propose in this work essentially depends on the ability to determine the constant $C_k$ which appears in the different corresponding approximation errors (\ref{estimation_error}), (\ref{erreur_ quadrature}) or (\ref{erreur_EDO}). \sa
Indeed, for each of these errors of approximation, the complexity of the constant $C_k$ will suggest appropriate investigations. In the current work we presented devoted to the finite elements method, we pointed out the important role of the \emph{a priori} estimates (\ref{estimation_a_priori}) and (\ref{estimation_a_priori_discret}), which are usually neglected, since they do not bring any asymptotic information/behavior when mesh size $h$ goes to zero.\sa
As a consequence, to render the probability distribution (\ref{F1})-(\ref{F7}) derived in Theorem \ref{Corollaire_Loi P_k1_k2_V2} operational, one will have to consider appropriate techniques which will make the determination (or at least the approximation) of the different constants which are involved in (\ref{F1})-(\ref{F7}) possible. This mainly owes to the contribution of the three estimates (\ref{estimation_a_priori}), (\ref{estimation_a_priori_discret}) and (\ref{estimation_error}) which leads to (\ref{Error_Estimate_Approx}), namely, $\hslash_{k_1}, \hslash_{k_2}$ and $h^*_{k_1,k_2}$. \\[0.4cm]
\textbf{\underline{Homages}:} The author wants to warmly dedicate this research to pay homage to the memory of Professors Andr\'e Avez and G\'erard Tronel, who broadly promoted the passion of research and teaching in mathematics.

\begin{thebibliography}{100}
%
\bibitem{Bramble_Hilbert} J. H. Bramble and S. R. Hilbert, Estimation of linear functionals on Sobolev spaces with application to Fourier  trnasforms and spline interpolation, \emph{SIAM J. Numer. Anal.}, {\bf 7}, pp. 112--124 (1970).
%
\bibitem{Brezis} H. Brezis, Analyse fonctionnelle - Th\'eorie et applications, Masson (1992).
%
\bibitem{ChaskaPDE} J. Chaskalovic, Mathematical and numerical methods for partial differential equations, Springer Verlag, (2013).
%
\bibitem{ChasAs18} J. Chaskalovic, F. Assous, A new probabilistic interpretation of Bramble-Hilbert lemma, \emph{Computational Methods in Applied Mathematics}, DOI: https://doi.org/10.1515/cmam-2018-0270 (2019).
%
\bibitem{Arxiv2} J. Chaskalovic, F. Assous, A new mixed functional-probabilistic approach for finite element accuracy, December 2018. arXiv:1803.09552 [math.NA]
%
\bibitem{arXiv_Wmp} J. Chaskalovic, F. Assous, Explicit $k-$dependenc for  $P_k$ finite elements in $W^{m,p}$ error estimates: application to probability distributions for accuracy analysis, Janvier 2019. arXiv:ZZZZ [math.NA]
%
\bibitem{Ciarlet} P.G. Ciarlet, Basic error estimates for elliptic problems, in Handbook of Numerical Analysis, Vol. II, Eds. P.G. Ciarlet and J. L. Lions, North Holland, (1991).
%
\bibitem{Crouzeix} M. Crouzeix et A. L. Mignot, Analyse num\'erique des \'equations diff\'ecentielles, Masson (1984).
%
\bibitem{Ern_Guermond} A. Ern, J. L. Guermond, Theory and practice of finite elements, Springer, (2004).
%
\bibitem{RaTho82} P.A. Raviart et J.M. Thomas, Introduction \`a l'analyse num\'erique des \'equations aux d\'eriv\'ees partielles, Masson (1982).

%
\bibitem{Stroud} A.H. Stroud, Numerical quadrature and solution of ordinary differential equations, Springer Verlag.(1974)
%
\bibitem{Toundykov_Avalos} D. Toundykov and G. Avalos, A uniform discrete inf-sup inequality for finite element hydro-elastic models, \emph{Evolution Equations and Control Theory}, {\bf 5(4)}, pp. 515-531 (2016).
%
\end{thebibliography}
\end{document}